\documentclass[12pt,a4paper]{article}
\usepackage{amsthm,amsopn,amscd,amssymb,amsfonts,amsmath,amstext}
\usepackage{fullpage}

\def\mytitle#1{\setcounter{equation}{0}                       %*
\setcounter{footnote}{0}                                      %*                                      
\begin{center}{\large\sc{#1}}\end{center}                     %*
\vspace{0.25cm}}                                              %*
\def\myname#1{\leftline{\bf #1}\vspace{-0.13cm}}              %*
\def\myplace#1#2{{\small\begin{flushleft}{\it #1}             %*
\quad{\tt#2}\end{flushleft}\vspace{0.1cm}}}                   %*
\renewenvironment{abstract}                                   %*
{\vspace{-0.5cm}\small                                        %*
\begin{center}                                                %*
 {\bf\normalsize Abstract}                                    %*
\end{center}}{}                                               %*
\newenvironment{keywords}                                     %*
    {\vspace*{3mm}                                            %*
    {\noindent{}\bf Keywords\/:}                              %*
        \nopagebreak\small}                                   %*
        {}                                                    %*
\newenvironment{ams}                                          %*
    {\vspace*{3mm}                                            %*
    {\noindent{}\bf                                           %*
    AMS 2000 subject classification\/:}                       %*
    \nopagebreak\small}                                       %*
    {}                                                        %*
\renewcommand{\thanks}[1]{                                    %*
    \renewcommand{\thefootnote}{~}                            %*
    \footnote{ #1 \vspace{0.09cm}}                            %*  
    \setcounter{footnote}{0}                                  %*
    \renewcommand{\thefootnote}{\arabic{footnote}}            %*
    }                                                         %*
\newenvironment{acknowledgement}                              %*
    {\section*{{\large Acknowledgement}}}                     %*

\newtheorem{theorem}{Theorem}[section]                        %*
\newtheorem{corollary}[theorem]{Corollary}                    %*
\newtheorem{lemma}[theorem]{Lemma}                            %*
\newtheorem{proposition}[theorem]{Proposition}                %*
\newtheorem{definition}[theorem]{Definition}                  %*
\newtheorem{conjecture}[theorem]{Conjecture}                  %*
\newtheorem{remark}[theorem]{Remark}

\begin{document}
%%%%%%%%%%%%%%%%%%%%%%%%%%%%%%%%%%%%%%%%%%%%%%%%%%%%%%%%%%%%%%%%%%%%%%%%%%%%%%%%%%%%%%%%%%%%%%%%%%
\mytitle{Measure convolution semigroups and non-infinitely divisible probability distributions.}
\myname{Aubrey Wulfsohn}
\myplace{Mathematics Institute, University of Warwick, Coventry, CV4 7AL, UK}{awu@maths.warwick.ac.uk}
\begin{keywords}
Probability distributions, Poisson processes, L\'evy processes, infinitely divisible, random variables, convolution semigroups, log-convex sequence, Stieltjes moment sequence, totally positive matrices, combinatorics,  random measures,  continuous products, Boolean convolution, statistical
mechanics.
\end{keywords}
\\
\\
\begin{ams}
60E10, 60E07, 05A05 (Primary) 46L51, 46N50, 82B10 (Secondary)
\end{ams}
\begin{abstract}
Let $\mu$ be a   probability measure (or corresponding random variable)
such that all moments $\mu_n$ exist. Knowledge of the moments is not
sufficient to determine infinite divisibility of the measure;  we show also that
infinitely divisible, and in particular lognormal, distributions lose infinitely
divisibilty when censored in certain ways even if all moments are arbitrarily
 close to those of the uncensored distribution. The moments of a composition
of k copies of $\mu$ are expressed as combinatorial compositions of the
$\mu_n$. We express the moments of the compositions in the context of
occupancy problems, arranging  n balls in k cells; the classical convolution
 is described by Maxwell-Boltzmann statistics and is multinomial. For certain
non-infinitely divisible measures with moments increasing  fast enough  the
indexing of a k-cell combinatorial composition is extended to indexing by
non-negative real t and we construct classical convolution measure semigroups
from amongst the t-indexed classes. We prove also that when a random variable
with infinitely divisible distribution is embedded in a L\'evy process $(Y_t)$
then the t-indexed Maxwell-Boltzmann is the law of $Y_t$. In order to get
moment-based multinomial compositions  indexed by a continuum we use random
measures and random distributions rather than random variables. An alternative
approach to  embeddability of a non-infinitely divisible $\mu$ is by
considering non-classical convolution measure semigroups; for example
embedding $\mu$ in a Boolean convolution measure semigroup and retaining the
multinomial character of the moments. Embedding  $\mu$ in an Urbanik
generalised convolution measure semigroup one loses the multinomial character
 of the moments.
\end{abstract}

\vskip 1000 pt

\section{Semigroups and non-infinitely divisible distributions}
\subsection{Moments of censored lognormal distributions}

The $(\alpha, \sigma)$-{\it{lognormal}} distribution $X_{\alpha,\sigma}$
has probability density \hfill\break
 ${ {1} \over { \sqrt {2\pi {\sigma}^{2}} y } }
{exp \lbrace 
 { { -(\log y - \alpha)^{2} }  \over { {2
{\sigma}^{2}} } } \rbrace }$, support
  ${\bf{R}}_{+}$
and is known to be infinitely divisible with moments $\mu_{n} = e^{n\alpha + {1 \over 2}n^{2}\sigma^{2}}$. A third parameter could describe shift of the support along ${\bf{R}}$. 
 There is no moment generating function since  the
  moment expansion of the 
characteristic function does not converge at all. The characteristic function 
is known (\cite{wul1}) and the 
density function for its  
L\'evy measure can be expressed as 
$\int_{0}^{\infty}{e^{-x \xi}U(\xi)d\xi}$  (\cite{wul2}), but $U$ cannot be explicitly 
determined.

For a real random variable  $X$  the {\it{spectrum}} $\sigma(X)$ of 
$X$ is taken to be  the set of all points of increase of its distribution function. It is well-known that
 the spectrum of an infinitely divisible
   distribution
   $\mu$ is  the closure of vector sum of $n$ copies of the spectrum of $\mu^{{1 \over n}}$, for every $n \in {\bf{N}}$.

For any lognormal distribution one can produce a discrete distribution with the same  moments which is not infinitely divisible. Indeed, it is shown in \cite{leip} that the discrete distribution with weights
 $a^{-n}e^{-{1 \over 2}n^{2}{\sigma}^{2}}$ at the points
  $a e^{n {\sigma}^{2}}, n \in {\bf{Z}}$
has the same moments as a lognormal distribution; the parameter $a$ is
 dependent on  $\alpha$. 
This distribution is not infinitely divisible;
 indeed, the spectrum of $\mu$ cannot be the  vector sum of  n copies of
    $\mu^{ 1 \over n}$, for any  $n > 1$.

We shall  generalise the following:
\begin{proposition}\label{p1} 
If $X$ is  a Poisson random variable  then 
 ${\bf{P}}[X \in (a,b)] > 0$ {\it{implies that}} 
 ${\bf{P}} [X \in (na,nb)] > 0$
 $ \lbrace$and also that\hfill\break
  $\sigma(X) \cap (a,b) \neq \emptyset$ implies
$\sigma(X) \cap (na,nb) \neq \emptyset \rbrace$ for all $n \geq 1$.
\end{proposition}

\begin{proof}
 Obvious since  ${\bf{P}}[ X = n] =
 e^{-\lambda} \lambda^{n} /n!$ for some $\lambda > 0$.
\end{proof}

\begin{definition}
The left-truncated lognormal is constructed  by truncating the mass over an interval $(0,b), b  > 0$ and
transferring the removed mass to the origin. The {\it{right truncated lognormal }} is constructed by truncating the interval $[c,+\infty)$ for some $c > 0$ and redistributing the removed mass over $[0,c)$. The gap-censored lognormal is constructed by cutting off the mass
over the interval $(a,b), 0<a<b$, and transferring that mass to the origin.
\end{definition}

\begin{proposition}
Right-truncated lognormal distributions are not infinitely divisible.
\end{proposition}
\begin{proof} They cannot be infinitely divisible as
  they have finite support ( \cite{luk} \S{8.4}). 
\end{proof}

\begin{definition} For a measure $\mu$ on ${\bf{R}}_{+}$ for which  all the moments are finite  we call  
$C =  limsup_{n} 
   {     {\mu_{n}^{2}} \over {{\mu_{n-1}\mu_{n+1}} }      }     $
  the critical ratio for the moment sequence of $\mu$.
\end{definition}
 In Proposition \ref{c41} the critical ratio  is shown to provide a criterion for indeterminacy of the moment sequence $(\mu_{n})$.
We  show  that 
censored lognormals can produce non-infinitely divisible distributions with the same critical ratio as a lognormal.

\begin{proposition} \label{p14}
For large $n$, by choosing small
  enough $b > 0$  the moment ${\tilde{m}}_n$ for the left-truncated lognormal 
distribution can
 be arbitrarily close to the untruncated $m_n$   and the critical
  value ${\tilde C}$ for the
   truncated moments will be arbitrarily close to $ C$.  The same is true for gap-censored lognormal distributions. 
\end{proposition}
\begin{proof}
 For the purposes of the Proposition there is no loss of generality
 in assuming $\alpha = 0, \sigma = 1$. Denote
 $\int_{x}^{\infty} {e^{-{1 \over 2}u^{2}}du} $ by
 $\Psi(x)$. Using \cite{ques}  one sees that 
 ${{{\tilde m}_{n}} \over m_{n} } =  \Psi(\log{b - n}) / \Psi (\log{b})$. Thus
${ {{\tilde {m}}_{n}} \over {m_{n}}  }$  is, for small enough $b$, arbitrarily 
close to 1 and $ { {\tilde{m}}_{n} } \over {{\tilde{m}}_{n \pm 1}}$ converges uniformly to $1$ as 
 $n \to +\infty$,  so 
 $\lim_{b \to 0} \limsup_{n}{ {{{\tilde{m}}_{n}}^{2}} \over {  {\tilde{m}}_{n-1} {\tilde{m}}_{n+1}} } $ =
  $\limsup_{n} { { {m_{n}}^{2} } \over { m_{n-1} m_{n+1} }   }$. The Proposition follows.
\end{proof}

\begin{definition}
 We call a Poisson distribution  with a random  intensity {\it{mixed Poisson}} and we   
denote the mixed Poisson with
intensity $X$ by ${\tilde{X}}$.
\end{definition}
 It is well-known that ${\tilde{X}}$ is
infinitely divisible if $X$ is infinitely
divisible.
\hfill\break

An infinitely divisible random variable $X$ on ${\bf{N}}_{0}$ has spectrum in ${\bf{N}}$ and probability
generating function $ \sum p_{n}z^{n} = e^{-\lambda (1 - Q(z))}$ where
 $p_{n} = {\bf{P}}[X = n]$ and $\lambda Q(z)$ is also a probability
generating function, written as $\sum_{0}^{\infty}{r_{k}z^{k}}$ with $r_{k} \geq 0$ for all $k$.
\begin{definition}\label{d17}
Katti's test \cite{kat} on a discrete distribution on
${\bf{N}}_{0}$ comprises  solving the equations $(n+1)p_{n+1} =
\sum_{k=0}^{j}{p_{j-k}r_{k},  j \in {\bf{N}}_{0}}$ for the given
$p_n$;  a negative $r_k$ implies that the distribution is not
infinitely divisible.
Conversely, the distribution is infinitely divisible if all the
$r_k$ are non-negative.
\end{definition}
\begin{proposition} \label{p13}
The left-truncated lognormal distribution is not always
 infinitely divisible; in particular it is not infinitely divisible when $\alpha =0,\hfill\break
 \sigma = 1, b =-1.4.$
\end{proposition}
\begin{proof}
  We do this by showing  that for the truncated distribution $X$ the discrete distribution ${\tilde{X}}$ is not infinitely divisible.    Consider a Poisson distribution with  a truncated lognormal 
intensity, 
so the probability generating function is determined by the $p_{k} = {{1} \over {k! \sqrt{2 \pi}}} \int_{\log{b}}^{\infty}{
e^{kx - e^{x} - {{x^{2}} \over {2}}}dx} $, and the truncated mass
${1 \over {\sqrt{2\pi}}} \int_{-\infty}^{-1.4}{e^{- {1 \over 2} {x^{2}}dx}}$ 
 is added to  $p_0$.
 Applying Katti's test  we obtain negative values of $r_{k}$ for 
sufficienly large  $k$. 

It was more economical for computation to  use  N copies of 
the truncated log-normal  setting
 $p_{k} =  {{N^{k}} \over {k! \sqrt{2 \pi}}}\int_{\log{b}}^{\infty}{
e^{kx - Ne^{x} - {{x^{2}} \over {2}}}dx} $. 
Independent computations on Maple and Mathematica, for
an accuracy of 16 digits and $N = 10$, produced negative $r_8$ and $r_9$.
\end{proof}

Katti's idea was extended, in \cite{sat}, to develop a general criterion for infinite divisibility of a
distribution on the half-line. We use this to prove:

\begin{theorem}
  The general left-truncated $X_{\alpha,\sigma}$
 is not infinitely divisible.
\end{theorem}
\begin{proof}
 A probability measure $\mu$ with support in ${\bf{R}}_{+}$ is not infinitely divisible
  if for some  $x > 0$
$$\int_{[0,x]}{y d\mu(y)} = \int_{(0,x]}{\mu([0,x-y])yd\nu(y)} + \gamma \mu([0,x])$$
for all $\gamma \geq 0$ and all   measures $\nu$ such that
$\int_{(0,+\infty)}{min(1,x)d\nu(x)} < +\infty$ \hfill\break
(\cite{sat} Theorem 51.1).
 Choosing $x < b$ then the left hand side
of the equation is zero whilst the right hand side is non-zero.
Indeed, $\mu[0,x]) = \mu([0,x-y])$ is the mass which had been
transferred from the truncated section to the origin, so $\gamma
\mu([0,x]) > 0$ and also
 $\int_{(0,x]}y d\nu(y) > 0$.
\end{proof}

\begin{remark}
Another method of proving non-infinite divisibility of distributions on the half-line is by determining  the  thickness of the tail. For example, truncate  the left half of the Gaussian distribution on ${\bf{R}}$; by adjusting the mass by moving  the truncated mass, say  to the origin, or by redistributing the left-hand-side mass over the right-hand-side.    The tail remains asymptotically the same as before adjustment  and so will be too thin for the distribution to be  an infintely divisible (see Theorem 26.1 of \cite{sat}).
\end{remark}
 \begin{definition}
We define a compound Poisson process to be of the form $X(t) =
\sum_{1}^{N(t)}{Y_{j}}$ where the $Y_{j}$ are i.i.d. non-negative
random variables representing the length of jumps taking place at
the epochs of a homogeneous Poisson process $N$, and $N(t)$
represents the number of epochs which have taken place by time
$t$.
In the literature the distribution has also been referred to as being
  completely random or as a mixed or generalised Poisson distribution. 
\end{definition}
The 
   characteristic function is known to
 be
 of the form \hfill\break
$\Phi_{X}(t) = exp\lbrace\int {(e^{i x  \theta} -1}) d \nu(x)\rbrace$  where
  its L{\'{e}}vy measure $\nu$  satisfies\hfill\break
$\int{min(1, x)d\nu(x)} < \infty$. The usual correction term in the expression 
is unnecessary since the integral converges without it. 

It is well-known that a distribution on ${\bf{R}}_+$ with an atom at $0$ is infinitely
divisible if and only if it is compound Poisson.

We shall usually interpret $X$ as the distribution for $X(1)$.

\begin{corollary}\label{c12}  If  the  jumps of a compound Poisson 
distribution $X$  are Poisson-distributed random variables
 its L{\'{e}}vy measure is discrete, and \hfill\break  
${\bf{P}}[ X \in (a,b)] > 0$  implies that
${\bf{P}}[X \in (na,nb)] > 0$
for all $n \geq 1$
\end{corollary}
\begin{proof} This is a direct consequence of Proposition \ref{p1}. 
\end{proof}
\begin{definition}
 For $X(t)$ be a compound Poisson process with L\'{e}vy measure
$\nu$ and $\epsilon
> 0$  we define $X_{\epsilon}$ to be the restriction of  $X$ to
jumps of length exceeding $\epsilon$.
\end{definition}
\begin{proposition} $\Phi_{X_{\epsilon}}(t)$  has  the form
 $exp \lbrace t \int_{\epsilon}^{\infty}{( e^{i \theta x} -1)
d\nu(x)} \rbrace$. 
\end{proposition}
\begin{proof} Partition $(\epsilon,+\infty)$ as $ \bigcup_{j=1}^{n}{B_{j}}$ 
where the $B_{j}$ are disjoint Borel sets.
  Define    
$X(B_{k},t)$ to be 
$ \sum_{ \tau \in \sigma(X)} {\chi_{B_{k}}{(X(\tau) - X(\tau-))}} $ , viz., the sum of the jumps of size in $B_{k}$ up to time $t$. For fixed $t$ the $X(B_{k},t)$
 are mutually 
independent and have characteristic functionals of the form
$exp \lbrace t \int_{B_{k}}{( e^{i \theta x} - 1) d\nu(x)} \rbrace$.
The Proposition follows.
\end{proof}
\begin{lemma}\label{l15}
 ${\bf{P}}[ \vert X - X_{\epsilon} \vert]> \eta >0$  converges to  0 as
$\epsilon \to 0+$.
\end{lemma}
\begin{proof} Since $\Phi_{X_{\epsilon}} \to \Phi_{X}$ as $\epsilon \to 0+$,
   by a standard argument $ \vert X - X_{\epsilon} \vert \to 0$ in probability.
 \end{proof}

\begin{definition} A random process $Z$ on (a product of complete separable metric spaces) 
 ${\cal{X}} \times {\cal{Y}}$ is called a {\it{marked point process}} if, for any Borel set $B$ in ${\cal{Y}}$ (the space of 'marks'), the number of points in the marginal process ${\hat{Z}}(B)$  is finite. 
\end{definition} 
 We shall show that $X_{\epsilon}$ is a marked point process, the  marks given by independent
 heights of the jumps. We need the following:

 If $\epsilon' < \epsilon$ then $\sigma(X_{\epsilon} ) \subset \sigma(X_{\epsilon'})$
 and the mass at any point of $\sigma(X_{\epsilon'})$ is larger than
the mass at $\sigma(X_{\epsilon})$. By Lemma \ref{l15}, for small enough $\epsilon > 0$ one has $[\sigma(X_{\epsilon} ) \cap (a,b)] \subset
[\sigma(X ) \cap (a, b)] $ and  the mass at any point of $\sigma(X_{\epsilon'})$
remains finite.

 \begin{proposition}
  Suppose  $X$ is compound Poisson with L\'{e}vy measure $\nu$ and choose $\epsilon > 0$. Let
 $N_\epsilon$ be a Poisson random variable with intensity $\nu([\epsilon,+\infty))$. Let
  $\nu_\epsilon$ denote the restriction of $\nu$ to $[\epsilon,+\infty)$. Marks are given
  by mutually independent heights $(h_{j}),\lbrace j = 1,\cdots, N_{\epsilon} \rbrace$, the
  jumps with distribution ${{\nu_{\epsilon}}
 \over {\nu(\epsilon,+ \infty})}$.
  Then $X_\epsilon$ is equal in law to $\sum_{j = 1}^{N_{\epsilon}}{h_{j}}$.
 \end{proposition}

\begin{lemma}\label{l18} 
 Let $(a,b)$ be a bounded interval in 
$(0,+\infty)$.
If \hfill\break
 ${\bf{P}}[X_{\epsilon} \in (a,b)] \geq 0$ then 
   for any
 integer 
$n \geq 1$ there is an arbitrarily small $\eta > 0$ such that also ${\bf{P}}[X_{\epsilon} \in (n(a - \eta),n(b + \eta)]
 \geq 0$.
\end{lemma}
\begin{proof}
 There is 
a finite set $J_{1}$ of integers  $j$ such that 
${\bf{P}}[\sum_{j \in J} {h_{j}} \in A] > 0$.
Fixing $n$, since the underlying process is Poisson, for every  
$h_{j}: j \in J_{1}$ there 
are also n different jumps with heights $h_{k}$ arbitrarily close  to 
$h_{j}$; we denote  
these heights by $h_{j}, j \in J_{\beta}, \beta = 2,3, ..,n$. Thus there will be an arbitrarily small
$\eta > 0$
such that  ${\bf{P}}[\sum_{j \in J_{1} \cup...\cup J_{n} } h_{j} \in
 (na - n\eta,nb +n\eta)] > 0$.
\end{proof}

\begin{theorem}
Suppose $0 < a <b $ and $X$ compound Poisson. Then ${\bf{P}}[X \in
(a,b)] > 0$ implies that ${\bf{P}}[X \in (na,nb)] > 0$ for any
integer $ n \geq 1$.
\end{theorem}

\begin{proof}
 Choose $ \eta > 0$. If ${\bf{P}}[X \in (a,b)] >0$ then  
${\bf{P}}[X_{\epsilon} \in (a-\eta,b+\eta)] >0$ for small enough $\epsilon$.
By Lemma \ref{p13} 
 one has
 ${\bf{P}}[ X_{\epsilon} \in (n(a - \eta), n(b + \eta)] > 0$. 
  Since  $X - X_\epsilon$ and $X_{\epsilon}$ are mutually independent, using
Corollary \ref{c12} and Lemma \ref{l18},
 $ {\bf{P}}[X \in (n(a- \eta),n(b+\eta))] \hfill\break
\geq {\bf{P}}[X_{\epsilon} \in (n(a-\eta),n(b +\eta)] \times {\bf{P}}
 [|X_{\epsilon} - X| < \eta]) > 0$.  Letting
  $\eta \to 0+$ the Theorem follows.
\end{proof}
\begin{corollary} Gap-censored compound Poisson distributions
are not \hfill\break
infinitely divisible.
\end{corollary}
\begin{proof}
Suppose there is a single gap between $a$ and $b$. If $2a < b$
then the interval $(a,2a)$ supports mass, a contradiction. If $2a
> b$ then the interval $(b,2a)$ supports mass,  if $2a = b$ then the interval
$(a/2,a)$ supports mass and again 
 there are contradictions. 
\end{proof}
\subsection{Combinatorial compositions and their convolution semigroups}
 
The convolution product of copies of a measure gives rise to a discrete convolution semigroup over 
${\bf{N}}_{0}$. We show that we cannot  embed  the measure in  a classical convolution semigroup over ${\bf{R}}_{+}$.

The n'th moment of a convolution product of copies of a probability measure is represented by ${\bf{E}}[(\sum_{i
=1}^{k} {X_{i}})^{n}]$ where the $X_{i}$ are independent identically distributed real random variables with
distribution $\mu$. We  admit non-classical interpretations of $(\sum_{i =1}^{k} {X_{i}})^{n}$.  

\begin{definition}
For  probability measures  such that all moments exist, $\mu \circ
\nu$  defines  the class of  measures having n'th moment $
\sum_{j=0}^{n}{ C^{n}_{j}   \mu_{j}\nu_{n-j}}$.  We call this the
binomial  composition of $\mu$ and $\nu$. It need not  define a
unique measure. We extend this idea to various combinatorial composition classes of $k$
copies of $\mu$ (not necessarily determining a unique measure), written say as $\mu^{\circ k}$. We denote the
classical convolution of $k$ copies of $\mu$ by $\mu^{k}$. We
express \hfill\break
 $(x_{1}+ .. +x_{k})^{n}$ as $ \sum_{ n_{1}+ .. +n_{k}=
n}{(n;n_{1},.. ,n_{k})x_{1}^{n_{1}}....x_{k}^{n_{k}}}$ where
$(n;n_{1},..,n_{k})$ denotes the number of ways of arranging $n$
objects, say balls or particles, into $k$  cells such that there
are $n_{j}$ balls in the j'th cell. The combinatorial composition
is called {\bf{multinomial}} if
 $$(n;n_{1},..,n_{q}) = { {n!} \over
{n_{1}! \cdots n_{q}!}}$$
\end{definition}
\begin{definition}
 A sequence of positive real numbers is 
called  Stieltjes if it is the moment sequence of some $\mu$ with support in ${\bf{R}}_+$, and Stieltjes
 determinate if this measure is  uniquely defined by the sequence.  
\end{definition}

\begin{theorem}
 Suppose  there is a combinatorial composition $\mu^{\circ k}, k \in {\bf{N}}$ of copies of
$\mu$ which can be extended to be indexed by $t \in{\bf{R}}_{+}$
in such a way that for  each $t$ the sequence $\mu^{\circ t}_{n}$
is a Stieltjes moment sequence, continuous at $t$ for each $n$.
Then there  is  a (classical) convolution measure semigroup
 $\Lambda(t)$ such that
$\Lambda(t)$ has the same moments as $\mu^{\circ t}$, though
$\Lambda(1) \neq \mu$ when $\mu$ is not infinitely divisible.
\end{theorem}
\begin{proof}
There is no loss of generality in assuming $t \in (0,1)$. For a
fixed $t$, since $\mu^{\circ t}_{n}$ is, by hypothesis, a
Stieltjes moment sequence, there exists a measure, which we denote
by $\lambda(t)$, with  moment sequence  $\mu^{\circ t}_{n}$. If
the moment-sequence is determinate then $\lambda(t) = \mu$. The
theorem follows from the following lemmas and the fact that if
$\mu$ is not infinitely divisible it is  not
stable and so not embeddable in the semi-group.
\end{proof}
\begin{lemma}
 Fix $t \in {\bf{R}}_{+}$ and denote the integral part of
$nt$ by $[nt]$. Given  a sequence of measures $(\lambda({1 \over
n})^{[nt]})$ there is a subsequence $(n^{(t)}_{k})$ of ${\bf{N}}$
such that the $ \lambda({1 \over n^{(t)}_{k}})^{[n_{k}t]}$
converges weakly to a measure, say $\Lambda(t)$, with the same
moments as $\mu^{\circ t}$ as $n \to \infty$.
\end{lemma}
\begin{proof}
This is standard, using Prokhorov's Theorem and the Chebyshev inequality.
\end{proof}
\begin{lemma}
 The $\Lambda(t), t \in {\bf{R}}_{+}$ form a semigroup.
\end{lemma}
\begin{proof}
In order to introduce sufficiently many  subsequences of sequences
of integers for all $s,t \in {\bf{R}}_{+}$ we start with $s$ and $t$ rational.
 Then $$\lim_{k}\lambda(\frac{1}{n_{k}})^{[n_{k}s]}*
\lim_{k}\lambda(\frac{1}{n_{k}})^{[n_{k}t]} *
(\delta_{(1/n_{k})})^{\rho}$$ converges to $\Lambda(s)
* \Lambda(t)$  for some subsequence $n_{k}$ of the sequences
$(n^{(s)}_{k})$ and $(n^{(t)}_{k})$. The factor $(\delta_{(1/n)})^\rho$, where $\rho$ can be $1$ or $0$, is there
to adjust the expression  whenever $nt \neq [nt]$. By hypothesis  $\mu^{\circ t}_{k}$ is continuous in a
neighbourhood of $t$ so ${\Lambda(t) * \Lambda(s) = \Lambda(s + t)}$ also  for $s,t \in {\bf{R}}_{+}$.
\end{proof}
\begin{remark}\label{r12}
If $\mu$ is not infinitely divisible and  the sequence $(\mu_{n})$ is determinate  then for any  combinatorial composition the $\mu^{\circ t}$ cannot be a Stieltjes moment sequence for all $t \in {\bf{R}}^{+}$. It was shown in \cite{bis} that the Maxwell-Boltzmann $(\mu^{\circ t}_{n})$, as defined the next section,
 form Stieltjes moment sequences if all the ${\mu_{n}^2} \over {\mu_{n-1}\mu_{n+1}}$ are all smaller than about $1 \over 6$.
\end{remark}

\section{Combinatorial convolution semigroups}
\subsection{Maxwell-Boltzmann statistics}
 We investigate combinatorial convolutions interpreting the $n$'th moment of combinatorial convolutions  $\mu^{\circ k}$ as  weighted combinations of the arrangements of  $n$ balls into $k$ cells so that $q = 0, 1, \cdots$ of these  are empty. To
determine the weights one needs some kind of statistics; there is no non-probabilistic formulation. Differences
arise if one distinguishes, or not, between the various particles labelled $x_{1}, x_{2}, \cdots$ and/or
distinguishes, or not, between the cells (say by selecting their order).

Consider a mechanical system  of $n$ balls and $k$ cells such that each ball is assigned to a cell. A {\it{state}} of the entire system is described by a random distribution of the particles in the cells. We
can distinguish between different phase-spaces, for example, with Maxwell-Boltzmann,  Bose-Einstein or Fermi-Dirac
statistics (see \S{3}). To apply the above to moments of measures we identify a cell with a measure $\mu$, and   $q$ balls in
$\mu$  is represented as $\mu_{q}$.

 There $k^{n}$ possible ways of distributing $n$ distinguishable balls in $k$ distinguishable cells. 
For the Maxwell-Boltzmann model, originally proposed to explain the distribution of sub-atomic particles into different energy states,  both the balls and the cells are
distinguishable and there are ${n! \over n_{1}!
\cdots n_{k}!}$  states; given cells numbered $1, \cdots, n$ and allowing multiple occupation and empty cells this is the number of  arrangements such that $n_{j}$ balls are in the $j'th$ cell.  Each of these states has  the
same probability of occurring i.e., $k^{-n}$. 
\begin{definition} 
$\mathfrak{B}^{n}_{k} = \sum_{j=0}^{k
-1}{ (-1)^{j} C^{k}_{j} (k-j)^{n}}$, zero for $n
> k > 0$. See \cite{Fel} II \S{11}. We shall call these Boltzmann numbers. 
\end{definition}
 There are $\mathfrak{B}^{n}_{k}$ 
distributions leaving no cells empty and the number of distributions leaving $q$ cells empty is
$C^{k}_{q}\mathfrak{B}^{n}_{k-q}$.  So also $C^{k}_{k-j}\mathfrak{B}^{n}_{j}$ is the number of distributions with only $j$ cells occupied. 

\hfill\break
{\bf{ Other characterisations of ${\mathfrak{B}}$.}}\hfill\break
 ${\mathfrak{B}}^{n}_{j} =
\sum_{n_{1}=1}^{n}\cdots \sum_{n_{j}=1}^{n} \chi_{(n_{1}+\cdots
+n_{j}=n)}\frac{n!}{n_{1}! \cdots n_{j}!}$, where $\chi$ signifies
the characteristic function of a set (cf. \cite{jor} \S{61}. \hfill\break
$\mathfrak{B}^{n}_{j} =
j!\mathfrak{S}^{n}_{j}$ where $\mathfrak{S}^{n}_{j}$ denotes the
Stirling subset number (i.e., the number of ways of partitioning $n$
objects in $j$ non-empty sets, identical to the Stirling number of
the second kind).\hfill\break
 $\mathfrak{B}^{n}_{j} = \Delta^{j}x^{n}_{|x = 0}$, where
$\Delta^{j}$ represents the $j$'th power of the finite difference
operator (cf. \cite{luk} \S{2.3}).
\begin{definition}
Given $\mu$, and  $k \in {\bf{N}}$, we define the Maxwell-Boltzmann $\mu^{\circ k}_{n}$   as
$$
\sum_{j=1}^{k}C^{k}_{j}\sum_{n_{1}=1}^{n} \cdots
\sum_{n_{j}=1}^{n} \chi_{(n_{1}+\cdots +n_{j}=n)} {n! \over n_{1}!
\cdots n_{j}!}\prod_{i=1}^{j} \mu_{n_i}
$$
For  $t$ not an integer  $C^{t}_{j}$  denotes the usual extension of $C^{n}_{j}$ for the  expansion of $(1 +
x)^{t}$. Since $C^{k}_{j} = 0$ when  $j > k$ we can define the Maxwell-Boltzmann $\mu^{\circ t}_{n}$ as
$$
\sum_{j=1}^{n}C^{t}_{j}\sum_{n_{1}=1}^{n} \cdots
\sum_{n_{j}=1}^{n} \chi_{(n_{1}+\cdots +n_{j}=n)}{n! \over n_{1}!\cdots n_{j}!}\prod_{i=1}^{j} \mu_{n_i}
$$
\end{definition}
\begin{proposition}\label{propn}
 The Maxwell-Boltzmann arrangement is multinomial.
\end{proposition}
\begin{proof}
 For $k \in {\bf{N}}$ the Maxwell-Boltzmann-$\mu^{\circ
k}_{n}$ is indeed $$ \sum_{ n_{1}+ .. +n_{k}= n}{ {n!} \over
{n_{1}! \cdots n_{k}!}}\prod_{i=1}^{j} \mu_{n_i}
$$
as the expression in the definition is only a re-arrangement of
the terms in this sum.
\end{proof}
 The moments are  well-defined,  involving finite sums of finite products of moments independent of $t$.
The  above definition can be confusing for computations. Seeing that the ${\mathfrak{S}}^{n}_{k}$ are well tabulated it is more convenient to compute  $\mu^{\circ t}_{n}$ as the polynomial 
 $$\sum_{j=1}^{n}{
 C^{t}_{j} \mathfrak{B}^{n}_{j} \sum_{n_{j}=1}^{n}{ \chi_{(n_{1}+\cdots +n_{j}=n)}} } \prod_{i = 1}^{j}{ \mu_{n_{i} }   }$$ 
We elaborate on this in \S{4}. Note that $\mathfrak{B}^{n}_{1} = 1$
and $\mathfrak{B}^{n}_{n} = n!$ for all $n$ and $\mathfrak{B}^{n}_{k} = 0$ for $k > n$. 

\hfill\break
 {\bf{Infinitely divisible distributions.}}

Since we are considering positive random
variables,  in order to simplify notation we shall occasionally use the (possibly
formal) moment generating function instead of the characteristic
function. Assuming
all moments of $\mu$ are finite, the cumulant generating function $\theta \mapsto \Psi(\theta)$, i.e., the log of the moment generating function,
can be expanded   as the infinite series \hfill\break
$\kappa_{1}\theta +
{{\kappa_{2}} \over {2!}}( \theta)^{2} + \cdots $ where
$\kappa_{i}$ is the $i$'th cumulant, viz. ${\Psi(0)^{(i)}} \over {i!}$ $ $. One 
compares this expansion with the infinite series obtained by expanding the characteristic function  into powers
of $\theta$; by using de Faa's formula for derivatives of a functional  one obtains expressions for the $\mu_{n}$ in as polynomials  in
 $(\kappa_{1}, \cdots \kappa_{n})$ (see \cite{luk} \S{2.4}).

\begin{definition}
Let $Y$ be a real random variable with infinitely
divisible distribution $\mu$.  We may assume that the
characteristic function of $Y$ is  of the form $ exp
\{ i \theta a + \int_{0}^{\infty}(e^{i\theta z} - 1) d\nu(z) \}$.
For $t \in {\bf{R}}_+$ let  $Y_{t}$ denote  the random variable
with characteristic function $ exp \lbrace t[ {i \theta a } +
\int_{0}^{\infty}{({e^{i \theta z} - 1) d\nu(z)}}] \rbrace$.
\end{definition} 
The characteristic function of an infinitely divisible $\mu$ is of
the form  $e^{\Psi(\theta)}$, so simplifying de Faa's expansion. One can show, using Taylor and McLaurin expansions  around $\theta = 0$,
 that 
$$\left.\frac{d^{n}}{
d\theta^{n}}\right|_{\theta = 0}(e^{\Psi(\theta)}) = n! \sum_{n_{1}
+\cdots + kn_{k} = n} (\Psi(0)^{(n_{1})} \cdots \Psi(0)^{(kn_{k})}))$$ 
giving an expression of  moments  in terms of cumulants.
The cumulant generating functional for $Y_{t}$ is thus $$t(\kappa_{1}\theta + {{\kappa_{2}} \over {2!}}( \theta)^{2} +
\cdots) $$
and the moments thus obtained  are
polynomials in $t$ with coefficients polynomials in the $\kappa_{i}$. 
\begin{theorem}\label{t25}
 For infinitely divisible $Y$ with distribution $\mu$ the 
classical  Maxwell-Boltzmann-$\mu^{\circ t}$  is the law
of the L\'evy process $Y_{t}$. 
\end{theorem}
\begin{proof}
The replacement of the symbol $k$ by $t$ in the Maxwell-Boltzmann
 $\mu^{\circ k}$  and the insertion of the multiplyer $t$ before the
 cumulant generating function of $Y$ to get the cumulant generating function of $Y_{t}$
are essentially trivial. Both $\mu^{\circ t}_{n}$ and
${\bf{E}}[(Y_{t})^{n}]$ are polynomials in $t$ vanishing at
integer values of $t$, and as such are identical.
\end{proof}
\hfill\break
{\bf{Illustration.}} \hfill\break
The simplest example is that  of a Poisson 
random variable  $Y$  with distribution $\mu$ and  mean $\lambda$. All the cumulants  equal $\lambda$. Using the difference-operator formulation of $\mathfrak{B}$ one
can show that $\mu_{n} =
\sum_{j =
1}^{n}{\mathfrak{B}^{n}_{j} \frac{\lambda^{j}} {j!}}$ and hence \hfill\break
${\bf{E}}[(Y_{t})^{n}] =
\sum_{j=1}^{n}{\mathfrak{B}^{n}_{j} \frac{(t\lambda)^{j}}{j!}}$. As  $ \mu_{1} = \lambda, \hskip 2pt  \mu_{2} = \lambda + \lambda^{2}, \hskip 2pt  \mu_{3} = \lambda + 3 \lambda^{2} + \lambda^{3}, \cdots$ \hfill\break
so ${\bf{E}}[(Y_{t})] = t\lambda, {\bf{E}}[(Y_{t})^{2}] = \lambda t +
 (\lambda t)^{2}, \cdots$
For  $t$ of the form ${1 \over k}$ one has, for example, \hfill\break
${\bf{E}}[Y_{1 \over 2}] = {1 \over 2}\mu_1, \hskip 2pt
 {\bf{E}}[(Y{_{1 \over 2}})^2] = {1 \over 2} \mu_2 - {1 \over 4} (\mu_1)^2,\hskip2pt
{\bf{E}}[(Y_{1 \over 2})^3] = {1 \over 2}\mu_3 - {3 \over 4} \mu_2\mu_1 + 
{3 \over 8} (\mu_{1})^{3}, \cdots$.
${\bf{E}}[Y_{{1 \over 3}}] = {1 \over 3} \mu_{1}, \hskip 2pt
 {\bf{E}}[(Y_{1 \over 3})^{2}] = {1 \over 3} \mu_{2} - {2 \over 9} (\mu_{1})^{2},\hskip3pt
{\bf{E}}[(Y_{1 \over 3})^{3}] =
 {1 \over 3} \mu_{3} - {3 \over 4} \mu_{2} \mu_{1} + {20 \over 27} (\mu_{1})^{3}, \cdots$, which agree with the corresponding combinatorial expressions for the Maxwell-Boltzmann-$\mu_{n}^{\circ {1 \over k}}$ and these could be used to construct the expectations for  $Y_t$ with $t \in {\bf{Q}}$

   \subsection{A continuum of cells}

Since $\mu$ is {\it{a fortiori}} not
stable it is not possible to embed $\mu$  in any of these
convolution semigroups or to interpret  this as  as an arrangement
of objects in   a continuum of cells (say  in a lattice with
spacing decreasing to zero).
For a chain $(X_{i})_{ i = 1,... ,k}$ of random variables ${\bf{E}}[(\sum_{i =1}^{k} {X_{i}})^{n}]$ is well-defined. When the $X_{i}$ are i.i.d. we have linked it to arrangements of $n$ elements in $k$ cells.
A  limiting  continuum of cells would be in the form
 ${\bf{E}}[(\sum_{s \in [0,t]}{X_{s}})^{n}]$ but this can have a meaning  only if at most a finite number of $X_s$ are non-zero. However by taking $X$ to be a random measure or random distribution the continuum of cells can have  a meaning. \hfill\break

\begin{definition}
A random measure is  taken  to be  a mapping associating a real random variable $X_A$ to each element $A$ of a
family of Borel sets in ${\bf{R}}_+$, or associating a real random variable $X(\phi)$ to test functions $\phi$
which are Borel on ${\bf{R}}_+$. A random (Schwartz) distribution associates random variables to test functions
$\theta$  which are infinitely differentiable with compact support.
\end{definition}
Following \cite{feld}, a random measure $X$ is said to be weakly  decomposable if, for every sequence of disjoint Borel subsets
$A_{1},A_{2}, ....$, the $X_{A_{i}}$ are mutually independent and  $\sum {X_{A_{i}}}$ converges almost everywhere
to $X_{\bigcup A{n}}$. The distribution of a weakly decomposable $X$ is interpreted as  the continuous product of
the distributions of its components It follows from \cite{feld} Theorem 2.1 that for a  weakly decomposable $X$
each $X_A$ is infinitely divisible and $\sum_{ \tau \in (s,t)}{X_{\lbrace \tau \rbrace}} = X_{(s,t)}$. \hfill\break

If
$X_{\lbrace t \rbrace}$ is zero for each $t$ then the characteristic functional is of the classical L\'evy-Khinchin
form.\hfill\break

A function $f$ is called decomposable if for every finite partition $A_{1}, \cdots ,A_{n}$ of $[0,1]$ into
measurable sets $f$ may be written as  $f_{1}+ \cdots +f_{n}$, where each $f_{j}$ is measurable with respect to
the field of measurable sets generated by the $A_{j}$. The conjecture in \cite{feld} that,  for the $\sigma$-field of Borel sets in $ [0,1]$,
 there are enough decomposable functions to generate all measurable functions  was verified in \cite{tsi}  by
showing that $t \mapsto X_{[0,t)}$ is a L\'evy process. The characteristic functional factors into a product of
independent Gaussian parts of various dimensions and a Poissonian part (see \cite{feld} Theorem 4.1). \hfill\break

What kind of processes, denoted say by $Y(t)$, can be constructed from random measures and distributions? One such process is $t \mapsto Y(t) = X_{[0,t)}$ as above.
For random distributions, if  $t \mapsto M(t,\theta)$ is integrable with respect to   $t$ and  $\theta \mapsto
e^{\int_{0}^{t}{M(s,\theta_{s}) ds}}$ is positive definite (to ensure that it is the Fourier transform of a
Schwartz distribution) then   $e^{\int_{0}^{t}{M(s,\theta_{s}) ds}}$ is the characteristic function of ${\sum_{ s \in [0,t)}}{X(\lbrace \theta_{s} \rbrace)}$.  A (Gaussian) white noise  corresponds to the characteristic functional with $M(t,\theta) = - {1 \over 2}
\theta^2$. Also a
Poisson-type  continuous product can be constructed based on a Fermi-Dirac system (cf.\cite{str}).

 These above  are {\it{Fock}}, also called {\it{linearisable}}, continuous products. The more interesting  non-linearisable
continuous products have been investigated in \cite{tsi} using  partitions 
other than  the Borel $\sigma$-fields.

\hfill\break
{\bf{Moments.}} To get moments of the continuous product  we could compute
  ${\bf{E}}[(X_{t})^{n}|X_{s}]$ or
 ${\bf{E}}[(X_{t})^{n} X_{s}(\phi)]$  where say  $X(\phi) = \sum_{k=0}^{m}a_{k}X^{k}$ or 
 $ X(\phi) = 
e^{izX}$. In the latter case for $t \in [0,s]$ we may write \hfill\break
 ${\bf{E}}[(X_{t})^{n}e^{izX_{1}}]
=  {\bf{E}}[(X_{t})^{n} e^{izX_{s}}] {\bf{E}}[e^{iz(X_{s} - X_{t})}]$.
When ${\bf{E}}[e^{izX_{s}}]$ is of the form  $e ^{s\psi(z)}$ one can calulate the moments by differentiation of the \hfill\break
 characteristic functional. However we already know by Theorem \ref{t25} that  when  $Y_{t}$ constructed from  $X$  is a  L\'evy process the moments 
 ${\bf{E}}[(Y_{t})^{n}|Y_{1}]$ are multinomial.

\section{Non-classical convolution measure semigroups}

\subsection{The generalized convolution}

A generalised convolution introduced by K.Urbanik, is a binary operation, indicated here by $\cdot$, on non-negative probability
distributions, associative commutative  and distributive for  convex sums. The translation mapping by $ x \in
{\bf{R}}_{+}$ is denoted by $T_{x}$. It is assumed also that there exists a sequence of real numbers $c_{n}$ such that
$T_{c_{n}}\delta_{1}^{\cdot n}$ converges weakly to a  measure $\sigma_{\kappa} \neq \delta_{0}$, called the {\it{characteristic measure}}. Several examples are given in \cite{urb1}.

The convolution is called regular if it admits  a (generalised) characteristic function such that there is a
bijective correspondence between probability measures  and  characteristic functions. The translation mapping by $ x \in {\bf{R}}_{+}$ is denoted by $T_{x}$. It is postulated that there exists a sequence of real numbers $c_{n}$ such that $T_{c_{n}}\delta_{1}^{\cdot n}$ converges weakly to a {\it{ characteristic measure}} $\sigma_{\kappa} \neq \delta_{0}$.   The  characteristic
function of $\sigma_{\kappa}$ is $\int_{0}^{+\infty}
{  exp^{-t^{\kappa} x} d \sigma_{\kappa}(x)}$ where $\kappa$ is
called the {\it characteristic index}. The semigroups defined by the generalised convolutions correspond to those for $\cdot$ -L\'evy processes as described in \cite{thu}.

 For the classical convolution
 $\sigma_{\kappa} = \delta_{0}$ and $\kappa = 1$ and  the characteristic function  is the Laplace transform of the measure.
\begin{proposition}
No generalised convolution, other than the classical convolution, is
 multinomial.
\end{proposition}
\begin{proof}
  For $\kappa = 1$ and $k > n$,
 $\mu^{\cdot k}_{n} =
\sum_{r=1}^{n}{(-1)^{n+r} C^{k}_{r} C^{k - r -1}_{n-r} (\mu_{n}^{\cdot r}})$, as in \cite{urb2}. This is
compatible with the multinomial combination. For $\kappa \neq 1$ the expansion is incompatible with the
multinomial expansion.
\end{proof}
\subsection{Non-commutative systems}

A non-commutative probability theory can be constructed as a non-commutative algebra ${\cal{A}}$ of observables with a normalised  positive linear state, say $\phi$; when  ${\cal{A}}$ is commutative to $x \in {\cal{A}}$ there corresponds a measure $\mu_{x}$
 such that $\mu_{x}(\phi) = \phi(x)$. Bose-Einstein and Fermi-Dirac  statistics relate to non-commutative quantum systems. Combinatorial convolution semigroups of measures derived  from quantum 'random variables'  relate to consideration of operators on Fock 
spaces. Eliminating the (Hamiltinian-based) dynamics and  fixing the number of particles,  so  ignoring creation and annihilation, one obtains a commutative subalgebra of ${\cal{A}}$ generated by the number-operators. For 
 Bose-Einstein and Fermi-Dirac systems the cells correspond to energy levels, 'balls' are indistinguishable and only
distinguishable arrangements
are considered so $(n;n_{1},..,n_{q}) = 1$. For Bosons there will be $C_{n}^{k+n-1}$ states.  For Fermi-Dirac statistics,
 by Pauli's exclusion principle, a cell will be either empty or
contain one ball;  there are  less balls than cells so  there will be $C^{k}_{n}$ states; in the computation of
$(\sum_{i =1}^{k} {X_{i}})^{n}$ one postulates $X_{i}X_{j} = - X_{j}X_{i}$ for $i \neq j$.   Convolutions of Bose-Einstein and Fermi-Dirac random variables have been   constructed in  \cite{vonw} by consideration of graded vector spaces and bi-algebras. The Fermi-Dirac convolution can represent a
 one-dimensional spin $1 \over 2$ Heisenberg ferromagnet (the ferro-magnetism as proposed by Heisenberg) sitting in a Fock  space, whereas the Bose-Einstein model has integer spin.  Quasi-particles, called magnons, carry the magnetic domain field in this model.
\begin{remark}  Seeing that a continuum ferromagnet  does not exist (see \cite{dub} and \cite{str}) there is no chance that the Fermi-Dirac 
  convolution product $\mu^{\circ n}$ can be extended to  a
 $\mu^{\circ t}$ as we  do for Maxwell-Boltzmann statistics.
\end{remark}
 
\subsection{The Boolean convolution}

The Boolean convolution product, related to a non-commutative probability theory, is  denoted by $\uplus$ has been dealt with in \cite{SpW}. A measure  $\mu$ can be  identified, by its moment sequence,  with a (normalised positive linear) states on ${\bf{C}}<X>$   and $\mu {\uplus} \nu$ is identified
with   states on ${\bf{C}}<X_{1},X_{2}>$, where 
 ${\bf{C}}<X_{1},X_{2},...,X_{n}>$ denotes the ring, with
complex coefficients, of polynomials in noncommuting indeterminates $X_{1},X_{2}..,X_{n}$. The restriction in \cite{SpW} that $\mu$ has compact support is unnecessary and the results are valid under the condition that all
moments exist (but then a moment sequence need not be determining).

  In \cite{Boz} the Boolean convolution semigroup ${\bf{R}}_{+} \ni t \mapsto \mu^{\uplus t}$ was shown to exist. The
semigroup is not the same as the classical convolution semigroup;  the characteristic function of $\mu$ is a
'self-energy' based on  the non-linear Cauchy transform. It can be expressed as $a + \int_{{\bf{R}}_{+}}{ {{1 +
xz} \over {z - x}}d\tau(z)}$ for $a \in {\bf{R}}$ and some measure   $\tau$ on ${\bf{R}}$, this form being
conducive to defining infinite divisibility of arbitrary distributions. Indeed it is shown  in \cite{SpW} that all $\mu$, as above, are 'infinitely divisible'.

The Boolean convolution is related to a lattice of {\it{non-crossing}} partitions; these are used also in quantum-type convolutions (cf. \S{2}),
\begin{proposition}
  The Boolean convolution is not a  bilinear relation; indeed it is highly non-linear and it is  not one of Urbanik's generalised convolutions. The Boolean convolution is binomial but not  multinomial. However the moments of  the $\mu^{\uplus k}$ satisfy multinomial decompositions and  $\mu_{n}^{\uplus k} = \mu_{n}^{k}$ for $n, k \in {\bf{N}}$.
\end{proposition}
\begin{proof}  It not  distributive for convex sums; for example,  ${1 \over 2}(\delta_{a} \uplus {1 \over 2} \delta_{b}) \uplus \delta_{c}$ is a convex sum of  $\delta_{a + b + c \pm {1 \over 2}\sqrt{(a-b)^{2} + c^{2}}}$ (see \cite{SpW} Example 3).  For  $n > 2$  the measure $
\mu^{\uplus t}$  does not have the same n'th moments as any element of $\mu^{\circ t}$. For example, for random
variables $X, Y, Z$ one has \hfill\break $(\mu \uplus \nu)(XYX) = \mu(X)\nu(Y)\mu(X)$  so  \hfill\break $(\mu
\uplus \nu)_{3} = m_{3} + 2m_{2}n_{1} + 2n_{2}m_{1} + (m_{1})^{2}n_{1} + (n_{1})^{2} m_{1} + n_{3}$ \hfill\break
as opposed to the classical $(\mu * \nu)_{3} =m_{3} + 3m_{2}n_{1} + 3m_{1}n_{2} + n_{3}$. The remainder of the Proposition follows as the $\mu^{\uplus k}, k\in {\bf{N}}$, are mutually commutative. 
\end{proof}

\section{Log-convex sequences, Stieltjes moment sequences and total positivity}

\begin{definition}
Sequences $(a_{n})$ of positive reals such that the  
  $ {a_{n}^{2}} \over { a_{n-1}a_{n+1}} $
 is  $\leq 1$, or  $\geq 1$, are called log-convex, or log-concave, respectively. We shall  call ${{a_{n}}^{2}} \over {{a_{n-1}a_{n+1}} }$
a moment ratio and denote it by $\theta_{n}$. We shall call the sequence $\theta$-log-convex if $\theta_{n} \leq \theta$ for all $n$. We call a $\theta$-log-convex  sequence strictly log-convex when $\theta < 1$. We similarly index log-concave sequences.
\end{definition}

We may assume, without loss of generality, that $a_{0} = 1$. The sequences 
of moment ratios were used by T.J.Stieltjes in \cite{sti} Chapter II. 
\begin{remark}\label{r42}
a).  It is shown in Theorem 51.3 of \cite{sat}, using Katti's test \ref{d17}, that a probability measure  on ${\bf{N}}$ is infinitely divisible  if $(p_{n})$ is log-convex.
\hfill\break
b).
Log-convex and log-concave sequences are known to be unimodal in the sense that they  are initially decreasing (increasing) with $n$ until one or two nodes $n_{0}$ and $n_{0} +1$ and then  increasing (decreasing) with $n$.
\end{remark}
 For $(a_{n}), (b_{n})$  log-convex  and log-concave  with moment ratios less and greater  than $\theta_{1}, \theta_{2}$ respectively,   the sequence $(a_{n}b_{n})$ will be log-convex if ${{\theta_{1}} {\theta_{2}}} \leq 1$.

The binomial sequences and the sequences 
$({\mathfrak{S}}^{n}_{k})_{k}$ of Stirling subset numbers are well-known to be log-concave.

\begin{proposition} \label{p42}
Stieltjes moment sequences are log-convex.
\end{proposition}
\begin{proof} 
It is obvious that
$\mu_{n}$ is  monotone non-decreasing. Since, as shown in  \cite{hash},   $\mu_{n+1} \over \mu_{n}$ is monotone non-increasing,
 $\theta_{n} = { {\mu_{n} \over {\mu_{n-1}}} / {{\mu_{n+1}} \over
  {\mu_{n}}}} \leq 1$.
\end{proof}

The following proposition for log-convexity can be adapted also to log-concave sequences.
\begin{proposition}\label{p43}
 Suppose the sequence $(a_{n})$ of positive real numbers has
  $\theta_{n} \leq \theta \leq 1$ for all $n$. Denote $log$ $a_{n}$ by $f_{n}$. \hfill\break
a) When $\theta < 1$ one can  extend the sequence ($f_n$)
 to a smooth continuous convex  function $F$ on the relevant interval of
 ${\bf{R}}_{+}$ with $F'' \geq log \sqrt{{1 \over \theta}}$,\hfill\break
 $F(t) \geq e^{ t^{2} log \sqrt{1 \over \theta}}  + O(t)$ and so $a_{n} > e^{n^{2} log {1 \over \sqrt{\theta}}}$. Here $F'(t)$ diverges as 
 $t\to \infty$.
\hfill\break
b) When  $sup$ $\theta_{n} = 1$ the extension of the $f_{n}$ to $F$ is such that \hfill\break 
 $\lim_{t \to \infty}{F'(t)} = \rho < \infty$;    $F(t)$ is asymptotically linear so $a_{n}$ is asymptotically  the order of  $ e^{\rho n} $. 
 \end{proposition}
\begin{proof}
 We construct $F$  by fitting parabolas touching  tangents, determined by
chords, at the points $(n,f(n))$. The 'acceleration' of $F$ at $n$ is \hfill\break
${
1 \over 2}(\log{a_{n+1}} + \log{a_{n-1}} - 2 \log{a_{n})} \geq -log \sqrt{\theta}$.   \end{proof}                                                                                                         
\begin{remark}
Case (b) cannot occur for  Stieltjes moment sequences  unless the moment sequence is determinate. If the characteristic function is analytic then                              $\lim_{n}sup{{{a_{n}}^{2} \over a_{n-1}a_{n+1}}} = 1$;  there is  large class of such $\mu$ having
 $\lim_{n}{  {a_{n}}^{2} \over a_{n-1}a_{n+1} } = 1$ 
but this does not hold  in general (see \cite{hash}).
\end{remark}
\begin{corollary}\label{c46}
If $(a_{n})$ is  strictly log-convex then for $ k < n$ one has \hfill\break
 ${ {a_{k} a_{n-k}} \over {a_{n}} } < 1$.  If $(a_{n})$ is  $\theta$-log-convex then for $ k < n$ one has 
 ${ {a_{k} a_{n-k}} \over {a_{n}} } \leq \theta^{k(n-k)}$.
\end{corollary}
\begin{proof}
From the  concavity of $F$ the gradient of $F$ is non-decreasing
 as n increases.
 Considering  $f_{0}, f_{k} f_{n-k}, f_{n}$ the first part of the Corollary follows. The second inequality
holds because the sequence $(\theta^{{n^{2} \over 2}} a_{n})$ will be log-convex.
\end{proof}
\begin{proposition}
 Fixing $n$, the sequence $({\mathfrak{B}}^{n}_{k})_{k}$ is log-concave 
and  $ k \mapsto \mathfrak{B}^{n}_{k}$ can be extended to a continuous convex function on
 $(0,n)$.
\end{proposition}
\begin{proof} If $j < n-1$, then by \cite{lieb}, 
 ${{(\mathfrak{S}^{n}_{j})^{2}} \over
 {\mathfrak{S}^{n}_{j-1}\mathfrak{S}^{n}_{j+1}}} \geq {{j \over{j-1}} {{n-j+1} \over {n-j}}}$  so \hfill\break
${{(\mathfrak{B}^{n}_{j})^2} \over {{\mathfrak{B}^{n}_{j-1} \mathfrak{B}^{n}_{j+1}}}} \geq {{j^{2} \over  {j^{2} -1}}} {{n - j +1} \over {n-j}} >1 $. The second part of the Proposition follows from Proposition \ref{p43}
\end{proof}
We note that the ${{(\mathfrak{B}^{n}_{j})^2} \over {{\mathfrak{B}^{n}_{j-1} \mathfrak{B}^{n}_{j+1}}}}$ will be  close to 1 except for small values of  $n$ so
the sequences $\mathfrak{B}^{n}_{j} \prod_{i=1}^{j}{\mu_{n_{i}}}$ will almost always  be log-convex.
\begin{definition}
 A sequence $(a_{i})$ is called  positive-definite if \hfill\break
$\sum_{i,j = 0}^{n}{ {{c_{i}c_{j}a_{i+j}} \geq 0}}$ for any  $n$ and any sequence of non-zero $(c_{i})$, strictly positive-definite 
when  the sum is greater than 0, semi-definite if the sum can be zero.  \hfill\break
The Hankel matrix $[a_{0} \cdots a_{2n}]$ is defined as the \small{$(n+1)
 \times (n+1)$} matrix having leading term $a_{0}$ and  $a_{i,j} = a_{i+j}$ for all
 $i, j \leq n$. Its determinant will be written as $ \Delta_{n}$ when it is clear what sequence is involved. The determinant of the Hankel matrix $[a_{p},\cdots,a_{p + 2n}]$ will be denoted likewise by $\Delta_{n,p}$.
\end{definition}
 A necessary and sufficient condition, going back to Stieltjes \cite{sti}, that a sequence $(a_{n}), a_{n} >0$ be a Stieltjes moment sequence  is that both \hfill\break
$n \mapsto a_{2n}$ and $n \mapsto a_{2n+1}$ are positive-definite functions for all $n$, and  so also if and only if both $[a_{0}, \cdots, a_{2n}]$ and $[a_{1}, \cdots a_{2n+1}]$ are positive definite for all $n$. 
\begin{definition}
Let $M$ be a square matrix  with coefficients say $M_{i,j}$. Let $\alpha = (\alpha_{1}, \cdots \alpha_{k}), \beta = (\beta{_1}, \cdots \beta_{k})$ where the $\alpha$ and $\beta$ are made up of non-decreasing sequences of integers. Denote by $M[\alpha|\beta]$ the minor obtained using rows and columns numbered $\alpha$ and $\beta$ respectively. When  $\alpha$ and $\beta$ are made up  of consecutive integers we call  $M[\alpha|\beta]$ a Fekete minor.\hfill\break
  The matrix is called positive-definite if $\sum_{i,j}{c_{i}c_{j}M_{i,j}}\geq 0$ or equivalently if all the principal minors are non-negative.
\hfill\break
 It is called totally positive if every minor is non-negative. We shall use the term strictly totally positive if  the above minors are  positive, and totally semi-definite if any of the minors are null.
\end{definition}
\begin{theorem} M. Fekete \cite{fek} Satz II.
A necessary and sufficient condition for a matrix be strictly totally positive is that all Fekete minors are strictly positive.
\end{theorem}

 \begin{theorem}
Let $(a_{n})$ be a sequence of positive reals. When $(a_{n})$ is a Stieltjes moment sequence then its associated   Hankel matrix
 is totally positive. If all  the Hankel matrices $ [a_{p}, \cdots,   a_{p + 2m}]$ are totally positive then $(a_{n})$ is  an indeterminate Stieltjes moment sequence.  
\end{theorem}
\begin{proof} 
  By definition, a totally positive matrix is also positive-definite. For $(a_{n})$ a Stietjes moment sequence both $n \to a_{2n}$ and $n \to
a_{2n+1}$ are positive-definite.
Any Fekete minor with leading term  $a_{2m}$ can be identified with a principal minor by considering  the sequence $n \to a_{2n}$ with leading term $a_{2m}$  and any Feteke minor with leading term $a_{2m+1}$ with a principal minor for the sequence
 $n \to a_{2n+1}$. Similarly mapping from
  principal minors for either $n \to a_{2n}$ or $n \to a_{2n+1}$ one can map to all the  Fekete minors. Therefore the Hankel matrix will be totally positive
by  Fekete's theorem.
 \end{proof}

\begin{theorem} \label{t4}
 a). A Stieltjes moment sequence is indeterminate if and only if both $n \mapsto a_{n}$ and $n \mapsto a_{n+1}$ are stricly positive definite \cite{ham1} (cf. Theorem 2.18 of \cite{sho}). \hfill\break
b). A Stieltjes moment sequence is indeterminate if and only if  both  the sequences  $({{\Delta_{n}} \over {\Delta_{2,n-1}}})$ and $({{\Delta_{1,n} \over \Delta_{3, n-1}}})$ have non-zero limits as $n \to \infty$ \cite{ham3} (cf.\cite{mer} Lemma 3).  \hfill\break
c). Let $(\mu_{n})$ be a Stieltjes moment sequence. Consider the  \small{ 2}-dimensional section  $D_{m} = \lbrace (x,y) : (x,y,\mu_{2},\cdots,\mu_{2m})\rbrace$ of  the cone in \small{(2n+1)}-dimensional space formed by  vectors $(\mu_{0},\cdots,\mu_{2m})$.
Let 
  $c_{1,m}$  be such that  $\Delta_{1,2m}$ vanishes if one replaces $\mu_{1}$  by $c_{1,m}$. The sequence $(c_{1,m})$ is non-decreasing, bounded by $\mu_{1}$ and we denote $\lim_{m \to \infty}{c_{1,m}}$ by $c_{1}$. When $c_{1} = \mu_{1}$ the sequence is determinate. By \cite{mer} Theorem 2, the  sequence is  indeterminate if and only if $\mu_{1} > c_{1}$ and $(c_{1},\mu_{1}] \subset D$ where $D = {\bigcap}_{m} {D_{m}}$ denotes  the limit  parabolic region as defined there.
\end{theorem}
\begin{proposition}
 A measure on ${\bf{R}}_{+}$ with moment generating function \hfill\break
 convergent in some neighbourhood of the origin is Stieltjes determinate. However   a Stieltjes moment sequence may determine a measure which has no well-defined moment generating function.
\end{proposition}
\begin{proof}
The first part is obvious since the generating function will define the measure
uniquely. Distribution functions $f$ with determining Stieltjes moment sequences
but  no well-defined moment generating function are easily found using 
G.H.Hardy's (1917, well before the Carleman) criterion for determinacy  viz., for some $q \geq 1$ and $\delta > 0$ one has
$\int_{0}^{\infty}{{f(u)}^{q} e^{\delta \sqrt{u}}} du < \infty$ (see. \cite{sho} I \S{6}). For
example, this holds  whenever
$f(u)$ is of the form $C e^{-u^{\alpha}}, \hfill\break
u >  0, {1 \over 2} \leq \alpha < 1$, noting that here the Hardy integrals  are   Laplace transforms. There is no moment generating function since $\int_{0}^{\infty}{ e^{su}  e^{- u^{\alpha}} du} = \infty$. 
\end{proof}
 Hardy's criterion can also be adapted for discrete distributions. Other suitable moment sequences could be constructed,  using Theorem \ref{t4}(c) and  manipulating Stieltjes' Hankel matrices.
\begin{remark}
a). The critical index $\delta$ does not  provide a  necessary condition for a log-conves  sequence to be Stieltjes or for a Stieltjes moment sequence to be  indeterminate.  Lognormal distributions, as in \S{1}, do not satisfy the $\delta$-condition if $\sigma^{2} \geq log {1 \over \sqrt{\delta}} > log 2$ but are 
indeterminate  no matter how large the moment ratio.
\hfill\break
b).  By Proposition \ref{p42} a moment sequence is always log-convex. The advantage of log-convexity is that one needs to compute the strict positivity or semi-definiteness of only one of the shifted sequences.
\hfill\break
c).  The authors of  \cite{crcs}, studying total-positivity for matrices, not necessarily Hankel,  get a slightly smaller value for $\delta$. They do not use a general sufficiency condition but compute inequalities using  selected critical minors and relating the critical index to a   rate of decrease of minors away from the diagonal.\hfill\break
 d). Can  the critical moment ratio be shown to be $1 \over 4$?
\end{remark}  
\begin{corollary} \label{c41}
The indeterminacy of a Stieltjes moment sequence depends only the asymptotic  
values of $\theta_{n}$ as $n \to \infty$.
\end{corollary}
\begin{proof} When one knows that the sequence is a moment sequence, for indeterminacy by
part (a) or  part (b)  one need  consider only the asymptototic values of the  $\theta_{n}$ as $n \to \infty$. 
\end{proof}
See Proposition \ref{p14} regarding our use of this Corollary.
 \begin{remark}
The moments of lognormal distributions, described in \S{1}, form totally positive Hankel matrices. Indeed  the matrices  are easily seen to have generalised Vandermonde determinants, well-known to be totally positive.  
\end{remark}
\begin{corollary}
If $(a_{n})$ is positive semi-definite every partial sum \hfill\break
$\sum_{i,j}{c_{i}c_{j}a_{i+j}}$ for $i,j \leq n$ can be positive for any fixed $n$.
\end{corollary}
\begin{remark}
  Sufficient conditions for positive definiteness of a matrix  involving the  off-diagonal $M[(i,j)|(i,j)]$ being small enough are known. For instance, it is not hard to prove that if $\sum_{i \neq j}{ \sqrt{M[(i,j)|(i,j)]}} \leq 1$ then the matrix is positive definite. This relates to the   decrease of off-diagonal minors used in \cite{crcs}.
In the case of Hankel matrices it is simpler to depend on having a fast enough rate of increase of the $\mu_{n}$.
\end{remark}

 As in \cite{boas} a  sequence of positive elements is a Stieltjes moment sequence if  the growth rate of the $a_{n}$ is fast enough. Indeed the growth rate
 \hfill\break
  $a_{n+1} \geq ((n+1)! a_{n})^{n+1}$ suggested there is  fast enough  for the sequence to be indeterminate. It is essentially sufficient that 
    the ${a_{n+1} \over a_{n}} \geq ({{1} \over {\delta}})^{n}$. 
\hfill\break

\hfill\break
{\bf{The Maxwell-Boltzmann convolution semigroup.}}
\hfill\break

It will henceforth be sufficient to  restrict $t$ to be in  $(0,1)$. The phenomena evident in  the expansions of $\mu^{t}_{n}$,  illustrated in \S{2} for Poisson $\mu$, persist for the Maxwell-Boltzmann $\mu^{\circ t}_{n}$ for a general $\mu$ and for this we will be  referring below  to 
the the alternative expansion in \S{2} 

\begin{lemma}\label{l44}
Being unimodal (see Remark \ref{r42}) the sequence $(\mathfrak{B}^{n}_{k})$, for fixed $n$, is at first increasing and eventally decreasing. It is known (see
\cite{jor}), that 
\hfill\break
${{\mathfrak{S}}^{n}_{k+1}} \over {{\mathfrak{S}^{n}_k}}$
$= {({{k+1}\over k})^{n}} {1 \over {k+1}}$ $\lbrace$ $
 {{1 - (k+1){({k \over {k+1}})}^{n} + \cdots}} \over {{1 - k({k-1 \over k})^{n}} + \cdots}$ $\rbrace \leq {{({{k+1}\over k})^{n}}  {1 \over {k+1}}}$; it follows that \hfill\break
${\mathfrak{B}^{n}_{k+1} \over 
{\mathfrak{B}}^{n}_{k}} \leq ({{k+1} \over k})^{n} {1 \over {(k+1)^{2}}}$.
 For
$k$ fixed, $\mathfrak{S}^{n}_{k}$ is  asymptotically $k^{n} \over k!$ so $\mathfrak{B}^{n}_{k}$ is asymptotically $k^{n}$. In particular, $2^{n-1} -1 \leq \mathfrak{S}^{n}_{2} \leq 2^{n-1}$ so $2^{n}-2 \leq \mathfrak{B}^{n}_{2} \leq 2^{n}$ for all $n$.
\end{lemma}

Expanding the series,
 
$\mu^{\circ t}_n = t \mu_n + {t(t-1) \over 2} 
\frac{1}{2^n -2}{\mathfrak{B}}^n_2(C^n_1{\mu_{n-1}\mu_1}  +
C^n_2\mu_{n-2}\mu_2 + \cdots$ \hfill\break
+ $C^n_{n-1}\mu_1\mu_{n-1})  +     
  C^t_3 {\mathfrak{B}}^n_3 (A_3^1 \mu_{n-3}\mu_2\mu_1 + \cdots)  + \cdots 
+ C^t_n n!(\mu_1)^n$.

The term accompaning ${\mathfrak{B}}^{n}_{2}$ accounts for  when there are $n-2$ empty cells. The coefficients
 ${{1} \over {2^{n} -2}}C^{n}_{j}$ are the probabilities of the various arrangements $\mu_{n-j} \mu_{j}$,
 $j =1,2,\cdots, n-1$. The 
$\mu_{n-1}\mu_{1}$ is put in the leading position as it will be  larger  than the other products of moments  for a reasonable rate of increase. In the  term corresponding to more non-empty cells the probabilities for the products of three moments are more complicated. We have not elaborated on these probabilities, denoted $A_{k}^{i}$, but inequalities  between the ${\mathfrak{B}}^{n}_{k}$, analogous to those used in Lemma 2 of \cite{bis}, can be obtained by using Corollary \ref{c46} and Lemma \ref{l44}.
\hfill\break

\begin{lemma} \label{l45}
Let $(\mu_{n})$ be a strictly log-convex Stieltjes moment sequence. For $t \in (0,1)$ each $\mu^{\circ t}_{n}$ can be written as a polynomial with leading term
$t \mu_{n}$, subsequent terms alternating in sign, because the terms of  $C^{t}_{n}$ has alternating signs, and  decreasing in modulus.  After any term in the sequence the absolute value of sum of the subsequent  terms is smaller than the absolute value of the preceeding term.
\end{lemma}
\begin{proof}
 The absolute decrease of the terms is proved using Lemma 2 of \cite{bis} with $\epsilon = 1$;  Corollary \ref{c46} gives
$\mu_{n-k}\mu_{k} <  \mu_{n}$  and this can be used to reduce products of k moments to products of k-1 moments, and so on. The last part holds generally for alternating series with terms decreasing in modulus.
\end{proof}

 \begin{corollary}\label{c420}
Given a $\theta$-log-convex moment sequence,  one has \hfill\break
$t \mu_{n} \geq \mu^{\circ t}_{n} > (1 - \theta) t\mu_{n}$ for all $n$.
\end{corollary}
\begin{proof}
 By Lemma \ref{l45} we know that $\mu^{\circ t}_{n} \leq t \mu_{n}$ and
 \hfill\break
 $\mu^{\circ t}_n \geq t \mu_n + \frac{t(t-1)}{2} {1 \over {2^n-2}}{\mathfrak{B}}^n_2
 (\sum_{k=1}^{n-1}{ C^n_k{\mu_{n-k}\mu_k}})$.  We use
${\mathfrak{B}}^{n}_{2} \geq 2^{n}-2$. From  Corollary \ref{c46} we have $\mu_{n-k}
\mu_{k} \leq \theta^{k(n-k)}$ so
 $\sum_{k=1}^{n-1}{C^{n}_{k}{\mu_{n-k}\mu_{k}}} \leq 2 \theta \mu_{n}$. The Corollary follows. 
\end{proof}

 \begin{proposition}\label{p420}
 For a $\theta$-log-convex Stieltjes moment sequence $(\mu_{n})$ with small
 enough $\theta$, or equivalently for a large enough rate of increase in the moments,
for any $t \in {\bf{R}}_{+}$  the Maxwell-Boltzmann $(\mu^{\circ t}_{n})$ is
a Stieltjes moment sequence.
\end{proposition} 
 \begin{proof}
It is sufficient  to show that for a large enough rate of increase of the moments one has ${({\mu^{\circ t}_{n}})^{2} \over \mu^{\circ t}_{n-1} \mu^{\circ t}_{n+1}} \leq \delta$  for all $n$. We know also that $\mu_{n+1} \mu_{n-1} \geq (\mu_{n})^{2}$,
the latter being the leading term in the denominator.
After cancelling  out the common factor $t^{2}$  the leading term,   in the denominator, will
  for fast enough increase of the   $\mu_{n}$  dwarf all other terms so the ${( {\mu^{\circ t}_{n})^{2}}} \over {\mu^{\circ t}_{n-1} \mu^{\circ t}_{n+1}}$ will be small enough to ensure that $(\mu^{\circ t}_{n})$ be an (indeterminate)  Stieltjes moment 
sequence.

Furthermore,  by Corollary \ref{c420}
 ${ ({\mu^{\circ t}_{n})^{2}} \over { \mu^{\circ t}_{n-1} \mu^{\circ t}_{n+1}} }
\leq { {\theta} \over {(1-\theta)^{2}} }$. As in \cite{bis} it is sufficient that   $ {{\theta} \over {(1-\theta)^{2}}} \leq \delta$, from which one can calculate the maximal $\theta$ that ensures that $(\mu^{\circ t}_{n})$ is Stieltjes.
\end{proof} 
\begin{conjecture}
  The maximal $\theta$ from the proceeding proposition is slightly less than ${1 \over 6}$. Because, as in Corollary \ref{c420},  the inequality has been  obtained   using only the first two terms of the expansion of $\mu^{\circ t}_{n}$   we conjecture that the estimate can be improved by consideration of further terms in the expansion of $\mu^{\circ t}_{n}$.
\end{conjecture}

\begin{acknowledgement}
The author wishes to thank Dr.R.Tribe for many helpful conversations and
Professor L.Bondesson for  suggesting
the use of mixed Poisson distributions in Proposition \ref{p13}.
\end{acknowledgement}

\end{document}